\documentclass[11pt]{amsart}
\usepackage{amssymb,amscd,amsmath,amsthm}
\usepackage[all]{xy}

\newtheorem {Theorem}                    {Theorem}[section]

\newtheorem {Proposition}[Theorem]       {Proposition}
\newtheorem {Claim}      [Theorem]       {Claim}

\newtheorem*{Proposition*}   {Proposition}

\theoremstyle{definition}
\newtheorem  {Definition} [Theorem]{Definition}

\theoremstyle{remark}
\newtheorem{Remark}[Theorem]{Remark}

\numberwithin{equation}{section}

\begin{document}
\title{Diffeological Smoothness In Hodge Theory}
\author{Jiayong Li}
\email{jiayong.li@utoronto.ca}

\begin{abstract}
On a compact, oriented, Riemannian manifold, the Hodge decomposition theorem associates a smooth primitive to any exact smooth form $\omega$. In this paper, we show that given a smooth family of exact smooth forms $\omega(t)$, the family of associated primitives is also a smooth family with respect to $t$.

\end{abstract}
		
\maketitle

\section{Introduction}

Let $M$ be an $n$ dimensional compact, oriented, Riemannian manifold $M$ with metric $g$. Let $\Omega^p(M)$ denote the space of smooth p forms on $M$, $*: \Omega^p(M) \rightarrow \Omega^{n-p}(M)$ the Hodge star operator, and $d: \Omega^p(M) \rightarrow \Omega^{p+1}(M)$ the exterior differential operator. Each $\Omega^p(M)$ is equipped with an $L^2$ inner product, $\langle \alpha, \beta \rangle := \int_M \alpha \wedge *\beta$, and $L^2$ norm, $\|\alpha\|_{L^2}^2 := \langle \alpha, \alpha \rangle$.

\

Recall that the co-differential $\delta: \Omega^p(M) \rightarrow \Omega^{p-1}(M)$ is $\delta := (-1)^{n(p+1)+1}*d*$, and the Laplace-Bertrami operator $\Delta: \Omega^p(M) \rightarrow \Omega^p(M)$ is $\Delta := \delta d + d\delta$. This operator is a natural generalization of the Laplace operator on Euclidean space. We refer the readers to Chapter 6 of \cite{warner} for a more detailed exposition.

\

The Hodge decomposition theorem states $$\Omega^p(M) = \Delta(\Omega^p) \oplus H^p = d(\Omega^{p-1}) \oplus \delta(\Omega^{p+1}) \oplus H^p,$$  where $H^p := \lbrace \omega \in \Omega^p: \Delta \omega = 0 \rbrace$. 

\

Applying $d$ on both sides of the equality $d(\Omega^{p-1}) \oplus \delta(\Omega^{p+1}) \oplus H^p = \Omega^p(M)$, we get $d: \delta(\Omega^{p+1}) \rightarrow d(\Omega^p)$. One can show that $d|_{\delta(\Omega^{p+1})}$ is in fact a bijection between $\delta(\Omega^{p+1})$ and $d(\Omega^p)$. Therefore, the space of p+1 exact forms, $d(\Omega^p)$, can be identified with $\delta(\Omega^{p+1})$ via $(d|_{\delta(\Omega^{p+1})})^{-1}$, so this operator gives a choice of primitive of any exact form. For simplicity, let us denote $(d|_{\delta(\Omega^{p+1})})^{-1}$ by $d^{-1}$ with domain and range understood. 

\begin{Definition}
We say that a family of forms $\{ \omega(t) \} \subset \Omega^p$ is \textbf{smooth} if, in local coordinates $x_1, \ldots, x_n$, the coefficents of $\omega(t)$ depend smoothly on $t, x$. Here $t \in \mathbb{R}^l$ is the parameter of the family.
\end{Definition}

\

The main goal of this paper is to show:

\

\begin{Theorem} \label{main}
Given a smooth family of exact forms $\omega(t)$, $d^{-1}(\omega(t))$ is also a smooth family.
\end{Theorem}

\begin{Remark}
The fact that $d^{-1}$ sends smooth family to smooth family exactly means that $d^{-1}$ is a \textbf{diffeologically smooth} map. For definition of diffeology, see the book \cite{pi} of Patrick Iglesias-Zemmour.
\end{Remark}

This result is well known to analysts, but we have not been able to find its explicit statement in the literature. It is used in Moser's method for a family of closed forms (Theorem 2 of \cite{moser}). Because of its importance for geometers, we find it worthwhile to spell out the details. It is a pleasure to thank William Goldman, Fr\'ed\'eric Rochon, Fabian Ziltener, and Adrian Nachman for sharing their insights and giving helpful references. Very useful comments were also made by Michael Bailey and Brian Lee. Lastly, the author would like to express his deepest gratitude to Yael Karshon who patiently supervised this undergraduate research project and offered guidance and encouragement.

\

The paper is organized as follows: In Section 2, we introduce some machinery from partial differential equations, such as $C^k$ spaces and Sobolev spaces. In Section 3, we prove Theorem \ref{main}, assuming the Proposition \ref{smooth2} (diffeological smoothness is the same as smoothness as map into $C^k$ spaces). In the appendix, we give a proof of Proposition \ref{smooth2}.

\

\section{Preliminaries}

\begin{Definition}
Let $U$ be an open set in $\mathbb{R}^n$, $f: U \rightarrow \mathbb{R}^m$, and $\alpha = (\alpha_1, \ldots, \alpha_n)$ an integer vector where each $\alpha_i \geq 0$. Then we define $|\alpha| := \sum \alpha_i$ and $\partial_x^{\alpha}f := \partial_{x_1}^{\alpha_1} \cdots \partial_{x_n}^{\alpha_n}f$, where the partial derivative is taken component-wise and should be understood in the weak sense. For the definition of weak derivative, see 5.2.1 of \cite{evans}.
\end{Definition}

\

\begin{Definition} \label{ck}
Let $U$ be an open and bounded set in $\mathbb{R}^n$, and $f: \mathbb{R}^n \rightarrow \mathbb{R}^m$ compactly supported in $U$. The $C^k$ norm of $f$ on $U$ is 

$$\|f\|_{C^k(U)} := \sum_{|\alpha| \leq k} \sup_{x \in U} |\partial_x^{\alpha}f|.$$  

More generally, the $C^k$ norm can be defined for forms on a compact manifold $M$ with charts $\lbrace U_i, \phi_i \rbrace$ and with partition of unity $\lbrace \lambda_i \rbrace$ (each $\lambda_i$ has compact support in $U_i$). Notice that if $\omega$ is a p-form, $\phi_{i*} (\lambda_i \omega)$ is just a vector valued function compactly supported in $U_i$. We define the $C^k$ norm of $\omega$ as $$\|\omega\|_{C^k(M)} := \sum_i \| \phi_{i*} (\lambda_i \omega) \|_{C^k(U_i)}.$$

Now we define the space of p-forms of type $C^k$ on $M$ as 
$$C^k(M, \Lambda^p) := \lbrace \omega:  \text{$\omega$ is a p-form and }\|\omega\|_{C^k(M)} < \infty \rbrace.$$

This is a Banach space containing the smooth p-forms $\Omega^p(M)$. 
\end{Definition}

\

The above norm is the sum of all the sup norms of the partial derivatives. Replacing the sup norm by the $L^2$ norm, we get the Sobolev norm:

\

\begin{Definition} \label{hs}
Let $U$ be an open and bounded set in $\mathbb{R}^n$, and $f: \mathbb{R}^n \rightarrow \mathbb{R}^m$ compactly supported in $U$. The $H^s$ norm (Sobolev s norm) of $f$ on $U$ is 

$$\|f\|_{H^s(U)}^2 := \sum_{|\alpha| \leq s} \int_U |\partial_x^{\alpha}f|^2.$$ 

More generally, the $H^s$ norm can be defined for forms on a compact manifold $M$ with charts $\lbrace U_i, \phi_i \rbrace$ and with partition of unity $\lbrace \lambda_i \rbrace$. We define the $H^s$ norm of $\omega$ as $$\|\omega\|_{H^s(M)} := \sum_i \| \phi_{i*} (\lambda_i \omega) \|_{H^s(U_i)}.$$

Now we define the space of p-forms of type $H^s$ on $M$ as 
$$H^s(M, \Lambda^p) := \lbrace \omega:  \text{$\omega$ is a p-form and }\|\omega\|_{H^s(M)} < \infty \rbrace.$$

This is a Banach space containing the smooth p-forms $\Omega^p(M)$. 
\end{Definition}

\

\begin{Remark}
Even though Definition \ref{ck} and Definition \ref{hs} depend on choices of charts and partition of unity, the topology induced by these norm do not, as shown in 1.3.4 of \cite{gilkey}. In particular, notice that the $H^0$ norm is just $L^2$ norm, which is equivalent to the (chart invariant version) $L^2$ norm defined in the beginning of the Introduction section. For a more detailed introduction to the $C^k$ and the $H^s$ norm, we refer the readers to Chapter 1 of \cite{gilkey} and Chapter 5 of \cite{evans}.
\end{Remark}

\

\section{Proof of Theorem}

By the Hodge decomposition theorem, $\Delta: (H^p)^{\bot} \rightarrow (H^p)^{\bot}$ is invertible. Let $G$, the Green operator, denote the inverse of $\Delta$. It can be easily shown that $d^{-1} = \delta G$ (use Proposition 6.10 of \cite{warner}). Therefore it suffices to prove that the Green operator carries a smooth family to a smooth family. 

\

It is useful to relate the notion of smooth family to the notion of smoothness in Banach spaces. (For differentiability and smoothness in Banach spaces, we refer the readers to Chapter 1 of \cite{lang}):

\begin{Theorem} \label{smooth1}
Let $\omega(t)$ be a family of p forms on a compact manifold $M$, parametrized by $t \in \mathbb{R}^l$. $\omega(t)$ is a smooth family if and only if for all k the map $\omega: \mathbb{R}^l \rightarrow C^k(M,\Lambda^p)$ is smooth.
\end{Theorem}

\

The above theorem says that the smoothness of $\omega(t)$ as a family is equivalent to smootheness of $\omega$ viewed as maps into $C^k$ spaces. The former notion is defined in terms of local coordinates, and in the latter notion, $C^k$ norms are defined in terms of local coordinates as well. Thus Theorem \ref{smooth1} follows from the following proposition, which we prove in the appendix.

\

\begin{Proposition} \label{smooth2}
Let $U$ be a bounded open set in $\mathbb{R}^n$, and $f: \mathbb{R}^l \times U \rightarrow \mathbb{R}^m$ be such that for each $t \in \mathbb{R}^l$, $f(t, \cdot)$ is compactly supported in $U$. Then,
\begin{align*}
f: \mathbb{R}^l \times U &\rightarrow \mathbb{R}^m \text{\quad is smooth} \\
(t,x) &\mapsto f(t,x)
\end{align*}
is equivalent to
\begin{align*}
\text{for all k, }\varphi: \mathbb{R}^l &\rightarrow C^k(U) \text{ is smooth} \\
t &\mapsto f(t, \cdot)
\end{align*}
\end{Proposition}

\

\begin{Remark}
This proposition holds when the domains of $t$ and $x$ are some infinite dimensional vector spaces and the smoothess is in certain generalized sense; see Theorem 3.12 in \cite{kriegl} for the exact statement. For it involves some technical functional analysis, we choose to prove the proposition in an elementary way rather than quoting this theorem.
\end{Remark}

\

We also need results on continuity of the Green operator on Sobolev spaces, and a relation between the $C^k$ norm and the $H^s$ norm. For the following propositions, see Lemma 1.3.5 (c) and (d) of \cite{gilkey} for functions, and a remark in section 1.3.5 of \cite{gilkey} for validity of these results for sections of vector bundles.

\begin{Proposition}
The Green operator $G: H^{s-2}(M, \Lambda^p) \rightarrow H^{s}(M, \Lambda^p)$ is continuous for all $s \geq 2$. 
\end{Proposition}

\begin{Proposition}
If $s > k + n/2$, then there exist positive constants $C_1, C_2$, such that, for all $\omega \in \Omega^p$, $C_1\|\omega\|_{C^k} \leq \|\omega\|_{H^s} \leq C_2 \|\omega\|_{C^s}$
\end{Proposition}

\begin{Remark} \label{bound}
In the above proposition, the inequality $C_1\|\omega\|_{C^k} \leq \|\omega\|_{H^s}$ follows from the Sobolev embedding theorem, and $\|\omega\|_{H^s} \leq C_2 \|\omega\|_{C^s}$ is obvious since the $L^2$ norm can be bounded by the sup norm. Combining the above two propositions, we see that $\| G(\omega) \|_{C^k} \lesssim \| G(\omega) \|_{H^s} \lesssim \| \omega \|_{H^{s-2}} \lesssim  \| \omega \|_{C^{s-2}}$, for $s > k + n/2$.
\end{Remark}

\

\begin{Remark} \label{partial}
Before proving Theorem \ref{main}, we first note that according to Proposition 3.5 of \cite{lang}, for $F$ a Banach space, the smoothness of a map $\varphi: \mathbb{R}^l \rightarrow F$ is equivalent to the existence and continuity of $\partial^{\alpha}\varphi: \mathbb{R}^l \rightarrow L^{|\alpha|}(\mathbb{R},F)$ for all multi-index $\alpha$, where $L^{|\alpha|}(\mathbb{R},F)$ denotes the space of multilinear continuous maps from $\mathbb{R}^{|\alpha|}$ to $F$. 

Moreover, we identify $\partial^{\alpha}\varphi: \mathbb{R}^l \rightarrow L^{|\alpha|}(\mathbb{R},F)$ with $\partial^{\alpha}\varphi(\cdot)(1)_1\cdots(1)_{|\alpha|}:\mathbb{R}^l \rightarrow F$. It is an easy exercise that this identification is norm preserving.
\end{Remark}

Now we are ready to prove Theorem \ref{main}:

\begin{proof}
Let $\omega(t)$ be a smooth family of forms. By Theorem \ref{smooth1}, for all $k$ the map $\omega: \mathbb{R}^l \rightarrow C^k(M,\Lambda^p)$ is smooth. In particular, the partial derivative $\partial_{t_i}\omega: \mathbb{R}^l \rightarrow C^k(M,\Lambda^p)$ exists and is continuous. We claim that, for all $k$ positive, $\partial_{t_i}G(\omega(t)):\mathbb{R}^l \rightarrow C^k(M,\Lambda^p)$ exists and is equal to $G(\partial_{t_i}\omega(t)):\mathbb{R}^l \rightarrow C^k(M,\Lambda^p)$, and it is continuous.

\

Indeed,
\begin{align*}
&\lim_{h \rightarrow 0} \frac{\| G(\omega(t+h\mathrm{e_i})) - G(\omega(t)) - G(\partial_{t_i}\omega(t))h \|_{C^k}}{h} \\
=& \lim_{h \rightarrow 0} \frac{\| G[\omega(t+h\mathrm{e_i}) - \omega(t) - \partial_{t_i}\omega(t)h] \|_{C^k}}{h} \\
\lesssim &  \lim_{h \rightarrow 0} \frac{\| \omega(t+h\mathrm{e_i}) - \omega(t) - \partial_{t_i}\omega(t)h \|_{C^{s-2}}}{h} \\
=& 0,
\end{align*}
where $\mathrm{e_i}$ is the i-th standard basis vector of $\mathbb{R}^l$. The first equality is by linearity of $G$, the inequality is by Remark \ref{bound} (for all $s$, $t$ such that $s > k + n/2$), and the last equality is by the definition of partial derivative,  $\partial_{t_i}\omega(t)$. Therefore by definition, for all $k$ positive, $\partial_{t_i}G(\omega(t))$ and $G(\partial_{t_i}\omega(t))$ are equal considered as maps from $\mathbb{R}^l$ to $C^k(M,\Lambda^p)$. Its continuity can be shown similarly: $$\lim_{t \rightarrow t_0} \| G(\partial_{t_i}\omega(t)) - G(\partial_{t_i}\omega(t_0)) \|_{C^k} \lesssim \lim_{t \rightarrow t_0} \| \partial_{t_i}\omega(t) - \partial_{t_i}\omega(t_0) \|_{C^{s-2}} = 0.$$
Since each $\partial_{t_i}\omega: \mathbb{R}^l \rightarrow C^k(M,\Lambda^p)$ is smooth, we can repeat this method and show that $\partial_t^{\alpha}G(\omega(t)) = G(\partial_t^{\alpha}\omega(t))$ and it is continuous. 
The rest follows from Theorem \ref{smooth1}.
\end{proof}

\

\section{Appendix: Proof of Proposition \ref{smooth2}}

In this section we shall prove that smoothness as a family is equivalent to smoothness in Banach spaces.

\begin{proof}[Proof of Proposition \ref{smooth2}]

We begin the proof by making an easy claim (without proof) about continuity:

\begin{Claim} \label{continuity}
If for each $t$, $f(t,x)$ is continuous in $x$ and $$\forall t_0, \ \lim_{t \rightarrow t_0} \sup_x |f(t,x)-f(t_0,x)| = 0,$$ then $f(t,x)$ is continuous in $(t,x)$.
\end{Claim}

\

We now show

\begin{Claim}
Assume $\varphi: \mathbb{R}^l \rightarrow C^k(U)$ is smooth for every $k$. Let $f(t,x) = \varphi(t)(x)$. Then all the partial derivatives of $f$ exist. In fact $\partial_t^{\alpha}\partial_x^{\beta}f(t,x) = \partial_x^{\beta}[\partial_t^{\alpha}\varphi(t)](x)$, and more generally, $$\partial_t^{\alpha_p}\partial_x^{\beta_p} \cdots \partial_t^{\alpha_1}\partial_x^{\beta_1}f(t,x) = \partial_x^{\beta_1+ \cdots +\beta_p}[\partial_t^{\alpha_1+ \cdots +\alpha_p}\varphi(t)](x).$$ Moreover, all the partial derivatives of $f$ are continuous in $(t,x)$.
\end{Claim}

\begin{proof}
We prove $\partial_t^{\alpha}\partial_x^{\beta}f(t,x) = \partial_x^{\beta}[\partial_t^{\alpha}\varphi(t)](x)$ by using induction on $|\alpha|$. 

\

We start with $|\alpha| = 1$ and choose $k \geq |\beta|$. By differentiability of $\varphi: \mathbb{R}^l \rightarrow C^k(U)$, $$\lim_{h \rightarrow 0} \frac{\| \varphi(t+h \mathrm{e_i}) - \varphi(t) - \partial_{t_i}\varphi(t)h \|_{C^k(U)}}{h} = 0.$$

\

Then it follows from definition of $C^k$ norm that
\begin{align*}
0 = 
&\lim_{h \rightarrow 0} \sup_{x \in U} \frac{|\partial_x^{\beta}[\varphi(t+h \mathrm{e_i})](x) - \partial_x^{\beta}[\varphi(t)](x) - \partial_x^{\beta}[\partial_{t_i}\varphi(t)h](x)|}{|h|} \\
= &\lim_{h \rightarrow 0} \sup_{x \in U} \frac{|\partial_x^{\beta}f(t+h\mathrm{e_i},x) - \partial_x^{\beta}f(t,x) - h\partial_x^{\beta}[\partial_{t_i}\varphi(t)](x)|}{h}.
\end{align*}
Thus for any $x \in U$, $$0 = \lim_{h \rightarrow 0} \frac{|\partial_x^{\beta}f(t+h\mathrm{e_i},x) - \partial_x^{\beta}f(t,x) - h\partial_x^{\beta}[\partial_{t_i}\varphi(t)](x)|}{h}.$$
By definition of $\partial_{t_i}$, $\partial_{t_i}\partial_x^{\beta}f(t,x)=\partial_x^{\beta}[\partial_{t_i}\varphi(t)](x)$. The rest of the induction is similar. Hence $\partial_t^{\alpha}\partial_x^{\beta}f(t,x) = \partial_x^{\beta}[\partial_t^{\alpha}\varphi(t)](x)$.

\

Choose $k \geq |\beta_1+\beta_2|$. Applying the above result twice,
\begin{align*}
&\partial_x^{\beta_2}\partial_t^{\alpha_1}\partial_x^{\beta_1}f(t,x) \\
=& \partial_x^{\beta_2}\partial_x^{\beta_1}[\partial_t^{\alpha_1}\varphi(t)](x) \\
=& \partial_t^{\alpha_1}\partial_x^{\beta_2}\partial_x^{\beta_1}f(t,x).
\end{align*}

We conclude that 
\begin{align*} &\partial_t^{\alpha_2}\partial_x^{\beta_2}\partial_t^{\alpha_1}\partial_x^{\beta_1}f(t,x) \\ 
=&\partial_t^{\alpha_2}\partial_t^{\alpha_1}\partial_x^{\beta_2}\partial_x^{\beta_1}f(t,x) \\
=&\partial_x^{\beta_1+\beta_2}[\partial_t^{\alpha_1+\alpha_2}\varphi(t)](x).
\end{align*}
Then $$\partial_t^{\alpha_p}\partial_x^{\beta_p} \cdots \partial_t^{\alpha_1}\partial_x^{\beta_1}f(t,x) = \partial_x^{\beta_1+ \cdots +\beta_p}[\partial_t^{\alpha_1+ \cdots +\alpha_p}\varphi(t)](x)$$ can be proved inductively.

\

To complete the proof of the claim, we need to show $\partial_x^{\beta}[\partial_t^{\alpha}\varphi(t)](x)$ is continuous in $(t,x)$. Choose $k \geq |\beta|$. It follows from the definition of smoothness of $\varphi: \mathbb{R}^l \rightarrow C^k(U)$ that $\partial_t^{\alpha}\varphi: \mathbb{R}^l \rightarrow C^k(U)$ is continuous, namely 
$$0 = \lim_{t \rightarrow t_0} \|\partial_t^{\alpha}\varphi(t)-\partial_t^{\alpha}\varphi(t_0) \|_{C^k(U)}.$$
Then $$0 = \lim_{t \rightarrow t_0} \sup_{x \in U} |\partial_x^{\beta}[\partial_t^{\alpha}\varphi(t)](x)-\partial_x^{\beta}[\partial_t^{\alpha}\varphi(t_0)](x)|.$$
By Claim \ref{continuity}, $\partial_x^{\beta}[\partial_t^{\alpha}\varphi(t)](x)$ is continuous in $(t,x)$.
\end{proof}

\

Now we show the other direction:

\begin{Claim}
Assume $f: \mathbb{R}^l \times U \rightarrow \mathbb{R}^m$ is smooth. Let $\varphi(t)(x) = f(t,x)$. Then $\partial_t^{\alpha}\varphi: \mathbb{R}^l \rightarrow C^k(U)$ exists. In fact, $\partial_t^{\alpha}\varphi(t) = \partial_t^{\alpha}f(t,\cdot)$ and $\partial_t^{\alpha}\varphi: \mathbb{R}^l \rightarrow C^k(U)$ is continuous for all $k$.
\end{Claim}

\begin{proof}
We shall prove this by induction on $|\alpha|$.

\

We start with $|\alpha|=1$ and fix $t \in \mathbb{R}^l$. By Taylor's Theorem, for each fixed $x \in U$, $$\partial_x^{\beta}f(t+h\mathrm{e_i},x) = \partial_x^{\beta}f(t,x) + \partial_{t_i}\partial_x^{\beta}f(t,x)h + \frac{1}{2}\partial_{t_i}^2\partial_x^{\beta}f(t+h'(\beta,x)\mathrm{e_i},x)h^2,$$
where $h'$ is some number between $0$ and $h$ and depends on $\beta$ and $x$. Then \begin{align*}
&\lim_{h \rightarrow 0} \sup_{x \in U} \frac{|\partial_x^{\beta}f(t+h\mathrm{e_i},x) - \partial_x^{\beta}f(t,x) - \partial_{t_i}\partial_x^{\beta}f(t,x)h|}{h} \\
= &\lim_{h \rightarrow 0} \sup_{x \in U} |\frac{1}{2}\partial_{t_i}^2\partial_x^{\beta}f(t+h'(\beta,x)\mathrm{e_i},x)|h \\
\leq &\lim_{h \rightarrow 0} Ch = 0,
\end{align*}
where the inequality is due to $f(t,\cdot)$ having compact support, and the fact that $h'$ lies in a closed interval around $0$ when we restrict $h$ to a closed interval, so $|\frac{1}{2}\partial_{t_i}^2\partial_x^{\beta}f(t+h'(\beta,x)\mathrm{e_i},x)|$ is bounded for each $t$.

\

Therefore 
\begin{align*}
& \lim_{h \rightarrow 0} \frac{\| \varphi(t+h\mathrm{e_i}) - \varphi(t) - \partial_{t_i}f(t,\cdot) \|_{C^k(U)}}{h} \\
=& \lim_{h \rightarrow 0} \sum_{|\beta| \leq k} \sup_{x \in U} \frac{|\partial_x^{\beta}f(t+h\mathrm{e_i},x) - \partial_x^{\beta}f(t,x) - \partial_{t_i}\partial_x^{\beta}f(t,x)h|}{h} = 0.
\end{align*}
Thus $\partial_{t_i}\varphi(t) = \partial_{t_i}f(t,\cdot)$. The rest of the induction is similar. This proves $\partial_t^{\alpha}\varphi(t) = \partial_t^{\alpha}f(t,\cdot)$.

\

To finish the proof of the claim, we need to show the continuity of $\partial_t^{\alpha}\varphi$, namely, $$\lim_{t \rightarrow t_0} \|\partial_t^{\alpha}\varphi(t)-\partial_t^{\alpha}\varphi(t_0) \|_{C^k(U)} = 0.$$
By the above result, it suffices to prove $$\lim_{t \rightarrow t_0} \sup_{x \in U} |\partial_x^{\beta}\partial_t^{\alpha}f(t,x)-\partial_x^{\beta}\partial_t^{\alpha}f(t_0,x)| = 0,$$ but this follows from mean-value theorem and boundedness of $\partial_t(\partial_x^{\beta}\partial_t^{\alpha}f)(t,x)$.
\end{proof}

Combining the above two claims proves the proposition.
\end{proof}
\
\
\

\end{document}